\documentclass[11pt]{amsart}
\usepackage{amssymb,amscd}

\newcommand{\cH}{{\mathcal H}}

\newcommand{\fb}{{\mathfrak b}}
\newcommand{\fg}{{\mathfrak g}}
\newcommand{\fh}{{\mathfrak h}}

\newcommand{\ft}{{\mathfrak t}}

\newcommand{\mc}{{\mathbb C}}

\newcommand{\mr}{{\mathbb R}}
\newcommand{\p}{{\p}}

\newcommand{\A}{{\alpha}}
\newcommand{\B}{{\beta}}

\newcommand{\Lie}{\operatorname{Lie}}

\def\rs{rationally smooth}
\def\sv{Schubert variety}
\def\svs{Schubert varieties}

\def\sv{Schubert variety}

\def\iff{if and only if }

\def\stb{said to be }
\def\st{such that }

\def\pp{Poincar\'e polynomial}

\def\p{\varphi}
\def\rs{rationally smooth}

\newtheorem{theorem}{Theorem}

\newtheorem{lemma}{Lemma}

\def\theckbibliography#1{\par\bigskip
\begin{center}
{\normalsize \bf References}
\end{center}
\par
\noindent\list
 {[\arabic{enumi}]}{\settowidth\labelwidth{[#1]}\leftmargin\labelwidth
 \advance\leftmargin\labelsep
 \usecounter{enumi}}

 \sloppy\clubpenalty4000\widowpenalty4000
 \sfcode`\.=1000\relax}

\begin{document}
\begin{center}

\bigskip
{\large\bf $B$-submodules of $\fg/\fb$ and\\ Smooth Schubert Varieties in
$G/B$}

\bigskip{\sc 
James B. Carrell}
\end{center}

\bigskip
{\tiny {\bf Abstract} Let $G$ be a semisimple linear algebraic group over $\mc$
without $G_2$ factors,
$B$ a Borel subgroup of $G$ and $T\subset B$ a maximal torus. 
The flag variety $G/B$ is a projective $G$-homogeneous variety whose tangent space at
the identity coset is isomorphic, as a $B$-module, to $\fg/\fb$,
where $\fg =\Lie(G)$ and $\fb =\Lie (B)$.
Recall that if $w$ is an element of the Weyl group $W$
of the pair $(G,T)$, the \sv\ $X(w)$  in
$G/B$ is by definition 
the closure of the Bruhat cell $BwB$.
In this note we prove that 
$X(w)$ is nonsingular \iff the following two conditions hold:
1) its \pp\  is palindromic and 
2) the tangent space $TE(X(w))$ to the set $T$-stable
curves in $X(w)$ through the identity is a $B$-submodule of $\fg/\fb$. 
This gives two  criteria in terms of the combinatorics of $W$
which are necessary and sufficient for $X(w)$ to be smooth:
$\sum_{x\le w} t^{\ell(x)}$ is palindromic, and 
every root of $(G,T)$ in the convex hull
of the set of negative roots whose reflection is less
than $w$ (in the Bruhat order on $W$) has the property that
its $T$-weight space (in $\fg/\fb$) is contained in $TE(X(w))$. 
owever, these conditions don't characterize the smooth \svs\ when $G$ 
has type $G_2$.}

\section{Introduction}

\svs\ in the flag variety of a linear
algebraic group $G$ were originally defined by Chevalley
in a famous unpublished paper \cite{CHEV}, where it was remarked
offhandedly that all \svs\ were probably smooth. This was an oversight,
of course, because it was already 
known that there exist singular \svs\
in Grassmannians.
The question of actually determining the smooth
\svs\ in an arbitrary flag variety 
has subsequently been treated in many places.
For example, when $G$ is of type $A$, Lakshmibai
and Seshadri \cite{LS} computed the tangent spaces of all
the \svs, and subsequently Deodhar \cite{DEOD} used this to show
that in type $A$, the \rs\ \svs\ are all smooth.
Dale Peterson showed this result holds 
whenever $G$ is simply laced (see \cite{CK} for a proof). 
In \cite{LSAND}, Lakshmibai and Sandhya showed that,
in type $A$, the smooth \svs\ are exactly the ones whose 
defining permutation avoids a certain pattern,
and, in \cite{BP}, Billey and Postnikov
extended pattern avoidance to
all classical $G$. 
Other results concerning  globally smooth \svs\ in $G/B$ 
in various settings include \cite{BILLEY}, \cite{BM}, 
\cite{CK}, \cite{KUM}, \cite{POLO}, \cite{RYAN}. 
The main result of this note 
extends Peterson's criterion for smoothness
to the non-simply laced setting by adding a condition 
which is vacuously satisfied in the simply laced setting.

\section{Preliminary Remarks}
Before stating the our result, we will review some elementary definitions 
and well known facts about \svs. 
Let $G$ denote a semisimple linear algebraic group over $\mc$ with a
fixed Borel subgroup $B$ and maximal torus
$T\subset B$, and let $\fg,\fb$ and $\ft$ denote their
respective Lie algebras. Recall the Cartan decomposition 
$$\fg=\ft \oplus \sum_{\A \in \Phi} \fg_\A$$
of $\fg$ into $T$-weight spaces,
where $\Phi$ is the root system of the pair $(G,T)$.
The set of positive roots consists of the roots 
corresponding to the $T$-module $\fb$. This set is denoted by  $\Phi^+$. 
One has $\Phi=\Phi^+\cup \Phi^-$, where $\Phi^-=-\Phi^+$.
Thus, the set of $T$-weights on $\fg/\fb$ is $\Phi^-$.
The flag variety $G/B$ of $G$ is a 
$G$-homogeneous projective variety, and it is well known
its tangent space $T_e(G/B)$ at the identity coset $e=B/B$ is
isomorphic with  $\fg/\fb$. Note, that
the identity coset in $G/B$ is fixed by $B$, so $T_e(G/B)$
is in fact a $B$-module, and, by considering the $B$-equivariant
projection $\pi:G \to G/B$, one obtains that
$T_e(G/B)$ and  $\fg/\fb$ are  isomorphic as $B$-modules.
We will henceforth make the identification
$$T_e(G/B)=\fg/\fb =\sum _{\A<0} \fg_\A.$$

Let $W=N_G(T)/T$ be the Weyl group of
$(G,T)$, and recall that $W$ is a finite reflection group 
generated by the reflections $r_\A$ (of $\fh$) through roots
$\A$. By the Bruhat decomposition  $G=BWB,$
$G/B$ is the union of the
$B$-orbits of the cosets $wB=n_wB$ as $w$ ranges over $W$,
where $n_w\in N_G(T)$ is a representative of $w$.
One of the basic results
of \cite{CHEV} says that the $B$-orbit  $Bw\subset G/B$ is
isomorphic to affine space $\mc^{\ell(w)}$,
where $\ell$ is the length function on $W$.
The Zariski closure $X(w)$ of 
the $B$-orbit $Bw$ is called the \sv\
associated to $w$. Thus $\dim X(w)=\ell(w)$.
 Another basic result of 
\cite{CHEV} is that the Bruhat order $\le$ on $W$, 
defined combinatorially in terms of
reflections, is compatible with the natural geometric 
order on the $B$-orbits: $x\le w$ \iff $X(x)\subseteq X(w)$.
Moreover, the $T$-fixed point set $X(w)^T$ is precisely
the set $\{x\in W\mid x\le w\}$.
It follows from these remarks
that the \pp\ of $X(w)$, with respect
to ordinary rational homolgy, has the well known expression
\begin{equation}\label{eqPP}
P(X(w),t)=\sum _{x\le w} t^{2\ell(x)}.
\end{equation} 

\section{Statements of Results}
The purpose of this paper is to give a simple 
constructive criterion which describes
which \svs\ with palindromic \pp\ in $G/B$ are smooth
that holds unless $G$ has a $G_2$ factor and,
in fact, is false in $G_2/B$. The condition that
the \pp\ of a \sv\ $X(w)$ is palindromic is equivalent to the 
more difficult to formulate condition 
that $X(w)$ is \rs (cf. \cite{CP}). A variety is \stb
\rs\ at a point if it satisfies local Poincar\'e 
duality at the point and globally \rs\ if it is
\rs\ at every point  \cite{KL}. We will state 
the smoothness critierion in three successive ways.
The first involves only the 
linear span $\Theta(w)\subset T_e(X(w))$ of the reduced tangent cone to
$X(w)$ at the identity element. 

\begin{theorem}\label{thm1} Suppose $G$ has no $G_2$ factors and
$X(w)$ is a \sv\ in $G/B$
whose \pp\ is palindromic (i.e. $X(w)$ is \rs).
Then $X(w)$ is smooth \iff $\dim \Theta(w)=\ell(w)$. 
\end{theorem}

In type $A$, it follows readily from a result
of Lakshmibai and Seshadri \cite{LS}
that $\Theta(w)=T_e(X(w))$. By combining a result
of the author \cite{C} and Polo \cite{POLO}, this equality also holds
in type $C$. Therefore, Theorem \ref{thm1} follows easily 
in types $A$ and $C$ from the Borel Fixed Point Theorem. 
Indeed, since $X(w)$ is $B$-stable
and the identity coset $e$ is fixed by $B$, 
$X(w)$ is smooth \iff it is smooth at $e$.

The second formulation of the smoothness criterion
requires that we bring 
in the tangent space $TE(X(w))$ to the  set $E(X(w))$ of 
$T$-curves to $X(w)$ at $e$ and discuss 
its relationship with $\Theta(w)$. 
Unless otherwise stated, proofs of the 
assertions here are in \cite{CP}.
Let
$$TE(X(w))=\sum_{C\in E(X(w))} T_e(C).$$ 
Each $T$-curve in $G/B$ is smooth, and hence if
$C\in E(X(w))$, $T_e(C)=\fg_\A$ for some $\A<0$.
One also has that $\dim TE(X(w))=|E(X(w))|$.
\begin{lemma} Let $\Phi(w)$ denote the set of
$\A<0$ \st $T$ has weight $\A$ on $T_e(C)$ for some $C\in E(X(w))$.
Then 
$$\Phi(w)=\{\A <0\mid  r_\A  \le w\}.$$
\end{lemma}

Thus
\begin{equation}\label{eqTE}TE(X(w))= \sum_{\A \in \Phi(w)} \fg_\A \subseteq \Theta(w).
\end{equation}
By Deodhar's inequality, $|E(X(w))|\ge \ell(w)$, 
so $\dim TE(X(w)) \ge \ell(w)$. 
It follows that $\dim \Theta(w)\ge \ell(w)$ with
equality \iff $TE(X(w))= \Theta(w)$. 
Finally, if $X(w)$ is \rs\ at $e$, then $\dim TE(X(w))=\ell(w)$.
However, knowing $\dim TE(X(w))=\ell(w)$ does not guarantee
that $X(w)$ is \rs\ at $e$.


Note that if $G$ is simply laced, then Theorem \ref{thm1} is 
vacuously true since a result of Dale Peterson's says that
every \rs\ \sv\ in $G/B$ is smooth (see \cite{CK} for a discussion and proof). 
The assumption that $\dim \Theta(w)=\ell(w)$ is 
automatically guaranteed when $X(w)$ is \rs,
since $TE(X(w))=\Theta(w)$ for all $w$ in the simply laced case.
In other words, Theorem \ref{thm1} gives the additional condition
under which Peterson's result holds
for the non $G_2$ setting.

The second formulation of Theorem \ref{thm1}
doesn't involve $\Theta(w)$.
Since $X(w)$ is $B$-stable and the identity coset $e$ is fixed by $B$,  
$T_e(X(w))$ and $\Theta(w)$ are $B$-stable submodules of $\fg/\fb$.
On the other hand, $TE(X(w))$ in general isn't $B$-stable. 
\begin{theorem} \label{main} Suppose $G$ has no factors of type $G_2$,
and assume  the \pp\ of the \sv\ $X(w)$ in $G/B$
is palindromic, i.e. $X(w)$ is \rs. 
Then $X(w)$ is smooth \iff $TE(X(w))$
is a $B$-submodule of $\fg/\fb$. 
\end{theorem}


The proof of this version uses the result that
$\Theta(w)$ is the $B$-module span of $TE(X(w)$
(see \cite[Theorem X]{C}).
Hence if $X(w)$ is \rs\ at $e$ and $TE(X(w))$ is a 
$B$-submodule of $\fg/\fb$, then clearly $\dim \Theta(w)=\ell(w)$.

The third formulation uses a description of the 
weights occuring in $\Theta(w)$. This
description involves the following convexity condition.
Note that here, the root system $\Phi$
is assumed to lie in the dual space $\ft^*$.
\begin{theorem}\label{thCONVEX} Let $w\in W$ be arbitrary, and 
suppose $\cH(w)$ is the convex hull 
of $\Phi(w)$ in $\ft^*$ over $\mr$. Then 
\begin{equation}\label{eqTHETA}
\Theta(w)=\sum_{\A \in \cH(w)\cap \Phi^-} \fg_\A.
\end{equation}
Thus $TE(X(w))=\Theta(w)$ \iff $\Phi(w)=\cH(w)\cap \Phi^-$.
Moreover, if $TE(X(w))$ is $B$-stable, then 
$TE(X(w))=\Theta(w)$.
\end{theorem}

\proof All one needs to do is quote
Theorem 3 and Corollary 1 of \cite{C}.\qed

\medskip
Our final formulation of the main result involves only $\Phi(w)$.
\begin{theorem} \label{thm4} Suppose $G$ has no factors of type $G_2$,
and assume  the \pp\ of the \sv\ $X(w)$ in $G/B$
is palindromic, i.e. $X(w)$ is \rs.
Then $X(w)$ is smooth \iff $\Phi(w)=\cH(w)\cap \Phi^-$.
\end{theorem}

Let us conclude this introduction with a few
remarks about the inclusions
$$TE(X(w))\subseteq\Theta(w) \subseteq T_e(X(w)).$$ 
In the $ADE$ setting, $TE(X(w))=\Theta(w)$. 
This is due to two facts: first,
every $T$-line in the reduced tangent cone to 
$X(w)$ at $e$ arises as $T_e(C)$ for a unique
$C\in E(X(w))$, and, second, every root is long
(hence  $\Phi(w)=\cH(w)\cap \Phi^-$).


\section{Proofs of Main Theorems}

It is clear from the above discussion that we only need to prove one 
version of the main result. Perhaps surprisingly, it turns out
that the easiest version to deal with is the $B$-module version
of Theorem \ref{main}.
We already noted that the main results hold in types $ACDE$,
so it suffices to check the $B$ and $F_4$ cases.
We will see below that type $B$ is easy to conclude from a result
of Billey. Thus the only sticking point is type $F_4$.
In order to get around this difficulty, we use the notion
of a stellar root system as introduced in \cite{BP}.

Let us begin by recalling what a stellar
root system is.
A reduced root system $\Phi$ distinct from $A_1$
and $A_2$ is called 
{\em stellar} if its Dynkin diagram is star shaped. That is, 
there exists a vertex which
is on every edge. Thus the stellar root systems are 
$B_2, C_2, G_2, A_3, B_3, C_3, D_4$. 
Let $\mr \Phi$ denote the real subspace of $\ft^*$ generated  
by $\Phi$. A subroot system of $\Phi$ is by definition a subset 
$\Psi$ of $\Phi$ of the form  $\Psi=\Phi \cap V$, 
where $V$ is a subspace of $\mr \Phi$.  A subroot system is a root system. Given $w\in
W$, the {\em inversion set}  of $w$ is the set
$I_\Phi(w)=\Phi^+ \cap w(\Phi^-)$. The inversion set $I_\Phi(w)$ uniquely determines
$w$, and for any root subsystem $\Psi=\Phi \cap V$, 
$I_\Phi(w)\cap V$ is the inversion set of a unique element 
$v \in W_\Psi$, the Weyl group of $\Psi$. The {\em flattening map}
$${fl}:W\to W_\Psi$$
is the assignment $w\to v$. That is, $fl(w)=v$. 

Let $\Psi$ be a subsystem of $\Phi$,
and let $G_\Psi$ (resp. $B_\Psi$) be the subgroup of
$G$ generated by $T$ and the root subgroups $U_\A$, 
where $\A\in \Psi$ (resp. the $U_\A$ with $\A \in \Psi^+$).
The following result of Billey and Postnikov classifies the 
\rs\ (resp. smooth) Schubert varieties for arbitrary 
$G$ in terms of stellar subsystems.
\begin{theorem}\label{thmSTELLAR}  A \sv\ $X(w)$ in $G/B$ is \rs\ (resp. smooth)
\iff for every stellar subsystem $\Psi$ of $\Phi$,
the \sv\  in $G_\Psi / B_\Psi$ corresponding to $fl(w)$ 
is also \rs\ (resp. smooth), where $(G_\Psi, B_\Psi)$
is the  unique pair of subgroups of $(G,B)$ determined by $\Psi$.
\end{theorem}

Let now prove Theorem \ref{main}.
Suppose $X(w)$ is a \rs\ \sv\ in $G/B$, and
let $\Psi$ be a stellar subsystem of $\Phi$. 
If $v\in W_\Psi$, let $Y(v)$ denote the corresponding
\sv\ in $G_\Psi / B_\Psi$. Put  $v=fl(w)$.  Then $Y(v)$ is \rs. 
By the discussion in \cite[Section 4]{BB}, 
$$Y(v)= X(w) \cap (G_\Psi / B_\Psi).$$ 
We claim this implies that $TE(Y(v))$
is $B_\Psi$-stable. To see this, suppose
$\A \in \Psi^-$ is a weight of $TE(Y(v))$,
and suppose that $\B\in \Psi^+$ is such that
$\A +\B$ is also a weight of the tangent space $T_e(Y(v))$. 
In particular, $\A +\B\in \Psi^-$, so
it follows that $\A +\B$ is a weight
in $T_e(X(w))$. By assumption, $\A +\B$ a weight of $TE(X(w))$,
since $TE(X(w))$ is  $B$-stable. Thus there
exists a $T$-invariant curve $C$ in $X(w)$ 
\st $C^T=\{e,r_{\A+\B}\}$ with weight $\A +\B$ at $e$.
Since $\A+\B \in \Psi^-$, $C\subset  G_\Psi / B_\Psi$ as well,  
so $C\subset Y(v)$. Consequently $\A +\B$
is a weight of $TE(Y(v))$, hence
$TE(Y(v))$ is $B _\Psi$-stable. Since $Y(v)$ is \rs\
and $\Psi$ is stellar, 
it suffices to verify that $Y(v)$ is smooth if 
$\Psi$ is of type $B_n$ for $n=2,3$ or $C_3$. This can be
checked directly, but it's more efficient to 
use the following lemma.

\begin{lemma} Theorem \ref{main} holds when $G$
is of type $B$ or $C$. 
\end{lemma}
\proof We have already verified this for type $C$.
Thus let $X(w)$ be a  \rs\ \sv\ in $G/B$
\st $TE(X(w))$ is a $B$-submodule of
$T_e(X(w))$, where $G$ is type $B$.
By the remarks in Section 3, 
it suffices to show $T_e(X(w))=TE(X(w)).$
since
$\dim TE(X(w))=\ell(w)$.
If $X(w)$ is singular, then one can apply Billey's pattern avoidance 
criterion \cite[Theorem 3]{BILLEY}.
Namely, $w$ contains the pattern $\bar{2}\bar{1}$ in  
signed permutation notation. By the argument on p.113 of \cite{BL}, 
it follows that there exists an $\A \in \Phi(w)$ 
and a $\B>0$  \st $\fg_{\A +\B} \subset T_e(X(w))$ 
for which the inequality
$r_{\A +\B}\le w$ fails. This says that
$\fg_{\A}\subset TE(X(w))$, while $\fg_{\A+\B}\not\subset TE(X(w))$,
contradicting the assumption that $TE(X(w))$ is $B$-stable.
Thus $X(w)$ must be smooth at $e$,
consequently smooth. \qed

\medskip Theorem \ref{main} now follows from 
Theorem \ref{thmSTELLAR} and the lemma.  \qed

\medskip
Let us next prove Theorem \ref{thCONVEX}.
To establish (\ref{eqTHETA}), we note the following
result \cite[Theorem 2]{C}: For any $x\le w$, let
$\Theta(w,x)$ denote the linear span of the reduced
tangent cone to $X(w)$ at $x$. Let $\cH(w,x)$
be the convex hull of $\Phi(w,x)=
\{ \A\in \Phi \mid x^{-1}(\A)<0, ~r_\A x \le w\}$.
Then 
$$\Theta(w,x)\subset \sum_{\A \in \cH(w,x)\cap \Phi} \fg_\A,$$
and any $\gamma \in \cH(w,x)\cap \Phi$ which isn't
a $T$-weight of $\Theta(w,x)$ 
has the form $\gamma =\B+\epsilon \mu,$
where $\B \in \Phi(w,x)$, $\mu >0$, $\epsilon \in \{1,2\}$
 and $x^{-1}(\mu)<0$. Thus, if $x=e$, $\gamma$ cannot exist.
\qed

\medskip
The fact that, in general, $\Theta(w)$ is the $B$-module span
of $TE(X(w))$ is proved explicitly in Theorem 3 of \cite{C}.

\section{Two Examples}
The first example shows that Theorem \ref{main} fails without
the $G_2$ hypothesis. That is, 
there exists a \rs\ but singular \sv\ $X(w)$ in $G_2/B$
for which $TE(X(w))$ is a $B$-submodule of $\fg_2/\fb$. Recall
that all \svs\ in $G/B$ are \rs\ when $G$ has rank two.

\medskip
\noindent
{\bf Example 1.}
Let $\A$ and $\B$ denote respectively the negatives of
long and short simple roots for $G_2$ corresponding to $B$, 
and let $r=r_\A$ and 
$s=r_\B$ be the corresponding
reflections. Let $w=srsrs$ and consider $X(w)$. 
Now $\ell(w)=5$ and it is not hard to see that
$$\Phi(w)=\{\A,\B,\A+\B,\A+2\B, \A+3\B\}.$$
Thus $TE(X(w))$ is indeed a $B$-submodule of $T_e(X(w))$.  
However, it is well known that $X(w)$ is singular:
for example, see \cite{KUM}. \qed

\medskip
\noindent
{\bf Example 2.} In this example,
we consider a singular \sv\ in the flag
variety $SO(5)/B$ of type $B_2$.
Let $\A$ and $\B$ denote respectively the 
negatives of the long and short simple roots
as in the previous example, and 
let $w=srs$. We claim 
$\Phi(w)=\{\A,\B,\A+2\B\}$, so $\Theta(w)$
is the $B$-module with weights
$\{\A, \B, \A + \B, \A+2\B\}.$ 
Thus $TE(X(w))\ne \Theta(w)$, hence, by Theorem \ref{main},
$X(w)$ is singular. Note that in the signed 
permutation notation for the elements of $W(B_2)$
(cf. \cite{BILLEY}), $w=\bar{2}\bar{1}$.
It is well known and easy to see
that $w$ is the unique element of $W(B_2)$
such that $X(w)$ is singular.
\qed   

It would be interesting to determine which $B$-submodules
of $\fg/\fb$ are tangent spaces at the identity to a smooth \sv\ in $G/B$.
This might lead to an efficient counting procedure for enumerating the
smooth \svs.

{\footnotesize
\begin{theckbibliography}{00}

\bibitem{BILLEY} S.\ Billey: {\em Pattern avoidance and rational smoothness of Schubert
varieties}, Adv. Math. {\bf 139} (1998), no. 1, 141--156.

\bibitem{BB} S.\ Billey and T.\ Braden:{\em Lower bounds for Kazhdan-Lusztig polynomials from patterns}, 
Transform. Groups {\bf 8} (2003), no. 4, 321--332.

\bibitem{BL} S.\ Billey and V.\ Lakshmibai: {\em Singular loci of Schubert varieties},
Progress in Mathematics, {\bf 182} BirkhŠuser Boston, Inc., Boston, MA, 2000. 

\bibitem{BM} S.\ Billey and S. A.\ Mitchell: {\em  Smooth and palindromic Schubert varieties in affine Grassmannians} 
 J. Algebraic Combin. {\bf 31} (2010), no. 2, 169--216. 

\bibitem{BP} S.\ Billey and A.\ Postnikov: {\em Smoothness of Schubert varieties 
via patterns in root subsystems.} Adv. in Appl. Math. 34 (2005), no. 3, 447--466.

\bibitem{CP} J. B.\ Carrell: {\em The Bruhat Graph of a Coxeter Group, a
Conjecture of Deodhar, and Rational Smoothness of Schubert Varieties,}
Proc.\ Symp.\ in Pure Math.\ {\bf 56}, No.\ 2, (1994), Part 1, 53-61.

\bibitem{C} J. B.\ Carrell: {\em  The span of the tangent cone of a Schubert variety},
Algebraic groups and Lie groups, 51--59,  Austral. Math. Soc. Lect. Ser., 9, Cambridge Univ.
Press, Cambridge, 1997. 

\bibitem{CK}  J. B.\ Carrell and J.\ Kuttler: {\em Singular points of $T$-varieties in $G/P$
and the Peterson map}, Invent. Math. {\bf 151}  (2003), 353--379.

\bibitem{CHEV} C.\ Chevalley: {\it Sur les decompositions cellulaires des 
espaces $G/B$}, Proc. Symp. in Pure Math. {\bf 56} (1994), Part I, 1-25.

\bibitem{DEOD} V. V.\ Deodhar: {\em
Local PoincarŽ duality and nonsingularity of Schubert varieties},
Comm. Algebra {\bf 13} (1985), no. 6, 1379--1388.


\bibitem{KL} D.\ Kazhdan and G. Lusztig: {\em Schubert varieties and Poincar\'Ž duality.}
Geometry of the Laplace operator (Proc. Sympos. Pure Math., Univ. Hawaii, Honolulu, Hawaii,
1979),  185--203, Proc. Sympos. Pure Math., {\bf 36}, Amer. Math. Soc., Providence,
R.I., 1980.

\bibitem{KUM} S.\ Kumar: {\em  Nil Hecke ring and singularity of Schubert
varieties}, Inventiones Math., {\bf 123} (1996), 471--506.

\bibitem{LSAND} V. \ Lakshmibai and B. \ Sandhya: {\em Criterion for smoothness of Schubert varieties in ${\rm SL}(n) / B$},
Proc. Indian Acad. Sci. Math. Sci. {\bf 100} (1990) 45--52. 

\bibitem{LS} V.\ Lakshmibai and C. S. \ Seshadri: {\em Singular locus of a Schubert variety},
Bull. Amer. Math. Soc. (N.S.) {\bf 11} (1984), no. 2, 363--366. 

\bibitem{POLO}  P.\ Polo: {\em On Zariski tangent spaces of Schubert varieties, and a proof
of a conjecture of Deodhar}, Indag. Math. (N.S.) {\bf 5} (1994), no. 4, 483--493. 

\bibitem{RYAN} 	K. M.\ Ryan: {\em On Schubert varieties in the flag manifold of ${\rm Sl}(n, C)$}, Math. Ann. 
{\bf 276} (1987) 205--224.

\end{theckbibliography}

\bigskip
\noindent
{\tiny
James B.\ Carrell\\
Department of Mathematics  \\
University of British Columbia \\
Vancouver, Canada V6T 1Z2 \\
carrell$@$math.ubc.c{a}}

\end{document}